\documentclass{article}
\usepackage{arxiv}

\usepackage[utf8]{inputenc} 
\usepackage[T1]{fontenc}    
\usepackage[colorlinks=true]{hyperref}
\usepackage{url}            
\usepackage{booktabs}       
\usepackage{amsfonts}       
\usepackage{nicefrac}       
\usepackage{microtype}      
\usepackage{lipsum}		
\usepackage{amsmath}
\usepackage{enumerate}
\usepackage[shortlabels]{enumitem}
\usepackage{tikz}
\usepackage{algorithm}
\usepackage{booktabs}
\usepackage[round]{natbib}
\usepackage{algpseudocode}
\usepackage{tikz}

\usepackage{booktabs}
\usepackage{siunitx}

\usepackage{graphicx}    

\newtheorem{definition}{Definition}
\newtheorem{theorem}{Theorem}
\newtheorem{lemma}{Lemma}

\newtheorem{proof}{Proof}
\def\QEDopen{{\setlength{\fboxsep}{0pt}\setlength{\fboxrule}{0.2pt}\fbox{\rule[0pt]{0pt}{1.3ex}\rule[0pt]{1.3ex}{0pt}}}}
\def\QED{\QEDopen}
\def\proof{\noindent{\bf Proof}: }
\def\endproof{\hspace*{\fill}~\QED\par\endtrivlist\unskip}
\newtheorem{remark}{Remark}
\DeclareSymbolFont{ugmL}{OMX}{mdugm}{m}{n}
\SetSymbolFont{ugmL}{bold}{OMX}{mdugm}{b}{n}
\DeclareMathAccent{\wideparen}{\mathord}{ugmL}{"F3}

\title{A Localized Method for Multicommodity Flow Problem}

\date{August 23, 2025}	

\author{
  Pengfei Liu\\
  Building No.1, Zhonguancun Road, Haidian District\\
  Beijing, China \\
  \texttt{pengfeil89@foxmail.com} \\
}

\begin{document}
\maketitle

\begin{abstract}
	This paper introduces a novel theoretical framework and a suite of highly efficient, parallelizable algorithms for solving the large-scale multicommodity flow (MCF) feasibility problem. We reframe the classical constraint-satisfaction problem as an equilibrium search. By defining vertex-specific height functions and edge-specific congestion functions, we establish a new, intuitive optimality condition: a flow is feasible if and only if it corresponds to a zero-stable pseudo-flow, where all potential differences across the network are resolved. This condition gives rise to an edge-separable convex optimization problem, whose structure is inherently suited for massive parallelization.

	Based on this formulation, we develop a family of Potential Difference Reduction (PDR) algorithms. Our primary method, provably convergent, solves an exact quadratic programming subproblem for each edge in parallel. To address scenarios with a very large number of commodities, we propose two computationally cheaper heuristics based on adaptive gradient descent. Extensive numerical experiments on well-known benchmarks demonstrate the framework's remarkable performance.   This work provides a powerful new approach for tackling large-scale MCF problems, while also identifying the formal analysis of the convergence rate as a promising direction for future research.
\end{abstract}


\section{Introduction}\label{sec:Intro}

The \emph{multicommodity flow problem} involves designing the flow for several different commodities through a common network with varying edge capacities.
We are given a directed graph $G(V,E)$ with $n$ vertices and $m$ edges, a capacity function $u : E \rightarrow \mathbb{Q}^+$, and $K$ commodities. Each commodity $k$ is defined by an origin-destination pair $(s_k, t_k)$ and a demand $d_k$. The objective is to find a flow assignment that satisfies the demand for each commodity without violating edge capacity constraints. The constraints can be summarized as follows:

\begin{equation}\label{eq1}
    \begin{aligned}
     &\sum_{k\in \mathcal{K}}f_{ij,k}\leq u_{ij}, \forall (i,j)\in E\\
     & f_{ij,k}\geq 0, \forall k\in \mathcal{K}, (i,j) \in E,\\
      &\sum_{j\in \delta^+(i)}f_{ij,k}-\sum_{j\in \delta^-(i)}f_{ji,k}=\left\{
      \begin{aligned}
      &d_k, \text{ if } i=s_k\\
      &-d_k, \text{ if } i=t_k\\
      &0,\text{ if } i\in V - \{s_k,t_k\}\\
      \end{aligned}
      \right.\\
      & \forall k\in \mathcal{K}, i \in V\\
    \end{aligned}
    \end{equation}
where \(\delta^+(i)=\{j|(i,j)\in E\}, \delta^-(i)=\{j|(j,i)\in E\}\) and  \(\mathcal{K} = \{1,2,\cdots,K\}\). In this paper, we assume $s_k\neq t_k$. The first expressions are capacity constraints. The second are non-negative constraints, and the last are flow conservation constraints.

Conventional centralized MCF algorithms require a central entity with complete network knowledge to compute optimal flows. In large-scale distributed networks, gathering this global information incurs significant communication overhead and delays \citep{farrugia2023solving}. Moreover, centralized approaches lack robustness due to a single point of failure and adapt poorly to dynamic network conditions. They are also unsuitable for scenarios where vertices must make autonomous decisions based on local information, as is common in ad-hoc networks or systems with dynamic membership \citep{awerbuch2007distributed}.

With the advancement of AI, the architectural paradigm of modern high-performance computing has shifted decisively towards massive parallelism \citep{robey2021parallel, d2021hardware}. Classical optimization algorithms designed for single-core processors often fail to exploit this potential. This technological shift creates a compelling need for algorithms whose computational structure mirrors the distributed nature of the network itself, making them both localized and inherently parallelizable.

This paper introduces a novel algorithmic framework for the multicommodity flow problem that directly addresses the limitations of existing methods by eliminating any reliance on global information. Departing from traditional methods, our approach simultaneously relaxes both capacity and flow conservation constraints. This complete relaxation transforms the problem into an equilibrium-finding problem governed by local interactions, guided by the physical analogy of flow moving from high to low potential.

To formalize this intuition, we introduce a new theoretical framework built upon local functions for congestion, nodal imbalance (height), and a ”driving force” we term potential difference. This leads to the definition of a stable pseudo-flow. Our central theoretical result is that a feasible solution to the
MCF problem exists if and only if a zero-stable pseudo-flow exists, which provides a complete and locally verifiable optimality condition. From this, we designed the \emph{Potential Difference Reduction}  (PDR) Algorithm, a simple iterative process that adjusts flows based on local rule, making the process inherently parallelizable.

The network notation introduced here is summarized in Table~\ref{Notation}. Further notation is introduced as needed.

\begin{table}[!htbp]
\centering
\caption{Basic Network Notation} \label{Notation}
\begin{tabular}{ll}
\hline
\(G(V,E)\)&~~~~~a directed graph \\
 \(V\)&~~~~~vertex (index) set\\
 \(E\)&~~~~~edge (index) set\\
 \(\mathcal{K}\)&~~~~~set of commodities, i.e., \(\mathcal{K} = \{1,\cdots,k,\cdots,K\}\)\\
 \(f_{ij,k}\)&~~~~~flow of commodity \(k\) on edge \((i,j)\), \(\boldsymbol{f} = (\cdots,f_{ij,k},\cdots)\)\\
  \(f_{ij}\)&~~~~~flow  on edge \((i,j)\), i.e., \(f_{ij} = \sum_{k\in \mathcal{K}}f_{ij,k}\)\\
 \((s_k, t_k, d_k)\)&~~~~~\(s_k\) and \(t_k\) are the origin and destination of commodity\\
 ~&~~~~~\(k\), and \(d_k\) is the demand\\
 \(u_{ij}\)&~~~~~the capacity of edge \((i,j)\)\\
 \(\delta^+(i)\)&~~~~~\(\{j|(i,j)\in E\}\)\\
 \(\delta^-(i)\)&~~~~~\(\{j|(j,i)\in E\}\)\\
 \(\delta_i\)&~~~~~the degree of vertex \(i\)\\
 \(\psi_{ij}\)&~~~~~congestion function of edge \((i,j)\)\\
  \(h_{ik}\)&~~~~~height function of vertex \(i\) for commodity \(k\)\\
  \(\Delta_{ik}\)&~~~~~the amount by which the \(kth\) flow  into the vertex \(i\) exceeds that out of the vertex\\
 \(r_{ij}\)&~~~~the unused capacity of edge \((i,j)\)\\
 \(m\)&~~~~ the number of edges of graph \(G(V,E)\)\\
 \(n\)&~~~~ the number of vertices of graph \(G(V,E)\)\\
 \(\phi_{ij,k}\)&~~~~ the potential difference across edge \((i,j)\) for commodity \(k\)\\
 \hline
\end{tabular}
\end{table}

\subsection{Our Contribution}

This paper presents a novel theoretical framework and a suite of highly efficient, parallelizable algorithms for solving the multicommodity flow feasibility problem. Our contributions are multi-faceted and can be summarized as follows:

\begin{enumerate}
    \item \textbf{A Novel Optimality Condition:}  By defining \textbf{vertex-specific height functions} and \textbf{edge-specific congestion functions}, we establish a simple and intuitive optimality condition, which reframes the  MCF problem from satisfying a set of hard constraints to finding an equilibrium state.

    \item \textbf{An Edge-Separable Convex Formulation:}  We formulate the stable pseudo-flow search as a convex optimization problem, yielding an inherent edge-separable structure, which makes the approach ideal for parallel processing on modern GPUs.

    \item \textbf{A Convergent, Decomposed Algorithm:} Based on the edge-separable formulation, we propose the \emph{Potential Difference Reduction Algorithm via Exact Subproblem Optimization}, which iteratively solves a small quadratic program for each edge in parallel.
    
    \item \textbf{Scalable Heuristics and Empirical Insights:}  Two computationally cheaper heuristic algorithms  are developed based on adaptive gradient descent, and an enhanced variant, which incorporates a \textbf{momentum term}, is specifically designed to accelerate convergence.

    \item \textbf{An Open Question for Theoretical Analysis:} While we have proven the convergence of the proposed methods, a formal analysis of their convergence rate remains an open problem. We posit that establishing a \emph{tight theoretical bound} on the rate is a non-trivial challenge, likely dependent on key topological and spectral properties of the underlying network graph.

\end{enumerate}

\section{Related Work}

Multicommodity flow (MCF) problems have attracted significant attention since the seminal works of  \cite{ford2015flows} and  \cite{hu1963multi}. For a comprehensive survey of the problem, see \cite{wang2018multicommodity1,wang2018multicommodity2,salimifard2022multicommodity}.

Many specialized solution methods based on linear programming have been proposed to exploit the block-angular structure inherent in MCF problems \citep{kennington1977solving,geoffrion1970primal, shetty1990parallel,frangioni1999bundle,goffin1997solving,moura2018branch, zenios1995smooth, kamath1995improved, kapoor1996speeding,cohen2021solving}. Conversely, studies by \cite{itai1978two} and \cite{ding2022two} have demonstrated that general linear programs can be transformed into two-commodity flow problems through nearly-linear time reductions.

Approximation algorithms for the MCF problem are also numerous \citep{leighton1999multicommodity, brun2017penalized, garg2007faster, leighton1995fast, schneur1998scaling, shahrokhi1990maximum,babonneau2006solving,chen2023high}. \cite{awerbuch1993simple, awerbuch1994improved} proposed an approximation algorithm that uses local-control techniques similar to the preflow-push algorithm of \cite{goldberg1987solving, goldberg1988new}. Unlike previous methods that find shortest or augmenting paths to push flow, their algorithm uses an "edge-balancing" technique that sends a commodity across an edge \((i,j)\) if the commodity is more queued at vertex \(i\) than at vertex \(j\). The push-relabel algorithm by \cite{goldberg1987solving, goldberg1988new}, which relaxes flow conservation constraints, is considered one of the most efficient algorithms for the maximum-flow problem. Recent improvements in this area compute the max-flow via a sequence of electric flows \citep{spielman2004nearly, daitch2008faster, christiano2011electrical, kelner2014almost, madry2016computing, axiotis2022faster, gao2022fully, chen2022maximum,van2023deterministic}. Building on this, \cite{van2023faster} developed a high-accuracy algorithm for the multicommodity flow problem using single-commodity methods.

Other approaches include using network equilibrium models to solve MCF problems \citep{liu2019combinatorial, bui2022decomposition, yu2022variable} and pseudo-flow methods for the max-flow problem \citep{hochbaum2008pseudoflow, chandran2009computational, fishbain2010competitive}. Flow decomposition has also been explored, often leveraging the problem's decomposable nature by commodity \citep{hartman2012split, haeupler2024low}.

Furthermore, research has explored parallel and distributed approaches for solving the MCF problem \citep{awerbuch1993simple, awerbuch1994improved, zhang2025solving, castro2000parallel, farrugia2023solving, ghaffari2015near,awerbuch2007distributed}. However, the effectiveness of parallelization can be limited by synchronization overhead and the inherently sequential components of some algorithms.

\section{Optimality Condition}
\label{sec:headings}

\subsection{Stable pseudo-flow}

Unlike other methods, our method \emph{does not enforce the capacity  and flow conservation constraints}.  The method, however, maintains a \emph{pseudo-flow}, which
is a function \(\boldsymbol{f}: K \times E \rightarrow \mathbb{R}^+\). Let \(\Delta_{ik}\) be the amount by which the
flow of commodity \(k\) into the vertex \(i\) exceeds the flow out of it, i.e.,
\begin{equation}\label{eq3}
\Delta_{ik}=\sum_{j\in \delta^-(i)}f_{ji,k} - \sum_{j\in \delta^+(i)}f_{ij,k} + \hat{\Delta}_{ik},\forall i \in V, k\in \mathcal{K},
\end{equation}
where \(\delta^+(i)=\{j|(i,j)\in E\}\), \(\delta^-(i)=\{j|(j,i)\in E\}\) and \(\hat{\Delta}_{ik}\) is defined as:
\begin{equation}\label{eq3h}
\hat{\Delta}_{ik}=\left\{
\begin{aligned}
&~~~0 ~~~~~~~~ if ~ i\in V - \{s_k,t_k\}\\ 
&~~~d_k  ~~~~~~~ if~i=s_k \\
&-d_k~~~~~~~ if~i=t_k.
\end{aligned}
\right.
\end{equation}

For each commodity \(k\)  at each vertex \(i\) , a height function \(h_{ik} = h_{ik}(\Delta_{ik})\) is introduced, 
where \(h_{ik}(\cdot)\) represents the relationship between \(\Delta_{ik}\) and the height of vertex \(i\) for commodity \(k\). Specifically, \(h_{ik}\) is defined as
\begin{equation}\label{eq2}
h_{ik}(\Delta_{ik}) = \frac{\Delta_{ik}}{\delta_{i}},
\end{equation}
where \(\delta_i \) is the degree of vertex \(i\).

That is, the height of vertex \(i\) for  commodity  \(k\)  is defined as the excess of inflow over outflow of commodity \(k\) at that vertex, normalized by the degree of vertex  \(i\).

Additionally, a \emph{congestion function} \(\psi_{ij} = \psi_{ij}(f_{ij})\) for each edge is introduced, where \(\psi_{ij}(\cdot)\) represents the relationship between the flow and the degree of congestion for edge \((i,j)\). Specifically, \(\psi_{ij}\) is defined as
\begin{equation}\label{eq41}
\psi_{ij}(f_{ij})=\left\{
\begin{aligned}
&0 ~~~~~~~~~~~~~~&~~~if~f_{ij} \leq u_{ij} \\
&f_{ij} - u_{ij}~~~~~ &~~~if~f_{ij} > u_{ij},
\end{aligned}
\right.
\end{equation}
where \(u_{ij}\) is the capacity of edge \((i,j)\) and \(f_{ij} = \sum_{k\in \mathcal{K}}f_{ij,k}\). If the flow on  edge \((i,j)\) is less than the capacity, the congestion of edge \((i,j)\) is zero. Otherwise, the congestion of edge \((i,j)\)  is the amount by which the flow  exceeds the capacity. 

The greater the \(h_{ik}\) is, the higher  the  vertex \(i\) for  commodity \(k\) is.  The greater the \(\psi_{ij}\) is, the more congested the edge \((i,j)\) is. 
As is often said, water finds its level. Intuitively, if  the height difference between vertex \(i\) and vertex \(j\) for commodity \(k\) is greater than the congestion of edge \((i,j)\), the flow \(f_{ij,k}\) should be increased;  if the height difference is less than the congestion, the flow \(f_{ij,k}\) should be decreased. In the remainder of this paper, , we  provide a rigorous proof to demonstrate that  the feasible solution for the multicommodity flow problem can be obtained by this intuitive idea.

First, we introduce the concepts of \emph{ potential difference} and \emph{stable pseudo-flow}.

\begin{definition}
\label{def0}
For every edge \((i,j)\) and commodity \(k\), the \emph{potential difference}    \(\phi_{ij,k}\) is defined  as  the height difference between vertex \(i\) and vertex \(j\), minus the congestion of edge \((i,j)\), i.e.,
\begin{equation}\label{eq4}
\phi_{ij,k} = h_{ik} - h_{jk} - \psi_{ij}.
\end{equation}
\end{definition}

\begin{definition}
\label{def1}
By using \{ \(\psi_{ij}, \forall (i,j) \in E\)\}  as the congestion function and  \{ \(h_{ik}, \forall i \in V, \forall k \in \mathcal{K}\)\} as the height function, a pseudo-flow \(\boldsymbol{f}\) is called a \emph{stable pseudo-flow} if it satisfies the following conditions:
\begin{enumerate}[label=(\roman*)]
\item for any used edge of commodity \(k\), the height difference between vertex \(i\) and vertex \(j\) for commodity \(k\) is equal to the congestion of edge \((i,j)\), i.e., the potential difference \(\phi_{ij,k}\) is zero;
\item  for any unused edge of commodity \(k\), the height difference between vertex \(i\) and vertex \(j\) is less than or equal to the congestion of edge \((i,j)\), i.e., the potential difference \(\phi_{ij,k}\) is less than or equal to zero;
\end{enumerate}
\end{definition}
where an edge \((i,j)\)  is called \emph{used} by commodity \(k\) if there exists \(s_k-t_k\) flow on edge \((i,j)\), otherwise it is called \emph{unused}.

\begin{definition}
\label{def2}
 A stable pseudo-flow \(\boldsymbol{f}\) is called  \emph{zero-stable} pseudo-flow if  \(\psi_{ij}=0\) and \(h_{ik} = 0\), for all \((i,j) \in E,  i \in V,  k \in \mathcal{K}\). Otherwise, it is called \emph{nonzero-stable} pseudo-flow.
\end{definition}

From the definitions above, a zero-stable pseudo-flow is a feasible flow that satisfies Expression~(\ref{eq1}). Therefore, we have the following theorem:
\begin{theorem}
\label{theo1}
Given \(\{(s_k,t_k,d_k): k \in \mathcal{K}\}\) and edge capacities \(\{u_{ij}: (i,j)\in E\}\), the feasible region of Expression~(\ref{eq1}) is not empty if and only if there exists a zero-stable pseudo-flow.
\end{theorem}

In fact, if  a nonzero-stable pseudo-flow exists, there is no feasible solution for Expression~(\ref{eq1}). 

\subsection{Basic Formulation}
Let \(f_{ij}\) be the sum of the flow of all commodities on edge \((i,j)\) and \(\Delta_{ik}\) defined as Expression~(\ref{eq3}). Define the following programming:
\begin{equation}\label{eq5}
\begin{aligned}
{\bf min}\quad &z = \sum_{(i,j) \in E}\int_0^{f_{ij}}\psi_{ij}(\omega)d\omega \\
&+ \sum_{i,k}\int_{0}^{\Delta_{ik}}h_{ik}(\omega)d\omega \\
{\bf s.t}\quad & f_{ij,k}\geq 0, \forall k\in \mathcal{K}, (i,j) \in E\\
\end{aligned}
\end{equation}
where \(\psi_{ij}\) is the congestion function and \(h_{ik}\)  the height function.

In the above programming, the objective function is the sum of the integrals of
the edge congestion functions and the integrals of the vertex height functions. 
It should be noted that there are only non-negative constraints present and no capacity constraint 
or flow conservation constraint. According to the definitions of the congestion function denoted by \(\psi\) and the height function denoted by \(h\), the following lemma can be derived:

\begin{lemma}
\label{lemma0}
The feasible region of Expression (\ref{eq1}) is  not empty if and only if the minimum value of the objective function of Expression (\ref{eq5}) is zero. Conversely, the feasible region of Expression (\ref{eq1}) is   empty if and only if the minimum value of the objective function of Expression (\ref{eq5}) is greater than zero.
\end{lemma}

\subsection{  Equivalence }
To demonstrate the equivalence between the stable pseudo-flow and the optimal solution of
Programming~(\ref{eq5}), it has to be shown that any flow pattern that solves Programming~(\ref{eq5})
satisfies the stable conditions. This equivalence is demonstrated  by proving that the Karush-Kuhn-Tucker conditions for Programming~(\ref{eq5}) are identical to the stable conditions.

\begin{lemma}
\label{lemma3}
Let \(\boldsymbol{f}^\ast\) be a solution of Programming~(\ref{eq5}).  \(\boldsymbol{f}^\ast\) is the optimal solution of Programming~(\ref{eq5}) if and only if it satisfies the  Karush-Kuhn-Tucker conditions of Programming~(\ref{eq5}).
\end{lemma}

\proof
Firstly, the objective function of Programming~(\ref{eq5}) is  convex. 
Secondly, the inequality constraints of  Programming~(\ref{eq5}) are continuously differentiable  concave functions. 
Therefore, Karush-Kuhn-Tucker conditions are necessary and sufficient for the optimality of Programming~(\ref{eq5}) (see \cite{boyd2004convex}).

\endproof

Since Programming~(\ref{eq5}) is a minimization problem with non-negativity constraints, the  Karush-Kuhn-Tucker conditions of such formulation are as follows:

\begin{equation}\label{eq6}
\begin{aligned}
&{\bf Stationarity}\\
& \quad -\frac{\partial z}{\partial f_{ij,k}} = -\mu_{ij,k},\forall k\in \mathcal{K}, (i,j) \in E\\
&{\bf Primal~feasibility}\\
& \quad -f_{ij,k}\leq 0, \forall k\in \mathcal{K}, (i,j) \in E\\
&{\bf Dual~feasibility}\\
& \quad \mu_{ij,k} \geq 0,\forall k\in \mathcal{K}, (i,j) \in E\\
&{\bf Complementary~slackness}\\
& \quad \mu_{ij,k}f_{ij,k} = 0,\forall k\in \mathcal{K}, (i,j) \in E.\\
\end{aligned}
\end{equation}
Obviously,
\begin{equation}
\label{cpgrad}
\begin{aligned}
\frac{\partial z}{\partial f_{ij,k}}  &= \psi_{ij}(f_{ij}) \frac{\partial f_{ij}}{\partial f_{ij,k}} \\
&+ h_{ik}(\Delta_{ik}) \frac{\partial \Delta_{ik}}{\partial f_{ij,k}} + h_{jk}(\Delta_{jk}) \frac{\partial \Delta_{jk}}{\partial f_{ij,k}}\\
	&= \psi_{ij}(f_{ij}) + (-h_{ik}(\Delta_{ik})) + (h_{jk}(\Delta_{jk})) \\
	&= \psi_{ij} + h_{jk} - h_{ik}.\\
\end{aligned}
\end{equation}
Substituting the  expression above into the Stationarity expression in KKT conditions,
\begin{equation*}
 h_{ik} - h_{jk} - \psi_{ij}=  -\mu_{ij,k} ,\forall k\in \mathcal{K}, (i,j) \in E.\\
\end{equation*}

For any used edge of commodity \(k\), i.e. \(f_{ij,k} > 0\), by complementary slackness \( \mu_{ij,k}f_{ij,k} = 0\) in KKT conditions, we have \( \mu_{ij,k} = 0\). Therefore,
\begin{equation*}
 h_{ik} -h_{jk} - \psi_{ij} = 0,\forall k\in \mathcal{K}, (i,j) \in E.\\
\end{equation*}
That is, the potential difference between vertex \(i\) and vertex \(j\) for commodity \(k\) is zero for any used edge of commodity \(k\).

For any unused edge of commodity \(k\), due to \(\mu_{ij,k} \geq 0\), we have 
\begin{equation*}
 h_{ik} -h_{jk} - \psi_{ij} = -\mu_{ij,k}\leq 0, \forall k\in \mathcal{K}, (i,j) \in E.\\
\end{equation*}

This direct interpretation of the Karush-Kuhn-Tucker conditions allows us to formulate the optimality conditions in terms of the potential difference \(\phi_{ij,k}\), which are equivalent to the complementary slackness conditions for this problem.

\begin{definition}[Stable Flow Optimality Conditions]
\label{def:sfo_conditions}
A pseudo-flow \(\boldsymbol{f}\) satisfies the optimality conditions for Programming~(\ref{eq5}) if and only if for every commodity \(k \in \mathcal{K}\) and every edge \((i,j) \in E\):
\begin{enumerate}[label=(\roman*)]
    \item if the edge is used ($f_{ij,k} > 0$), the potential difference is zero:
    \[ \phi_{ij,k} = h_{ik} - h_{jk} - \psi_{ij} = 0; \]
    \item if the edge is unused ($f_{ij,k} = 0$), the potential difference is non-positive:
    \[ \phi_{ij,k} = h_{ik} - h_{jk} - \psi_{ij} \leq 0. \]
\end{enumerate}
\end{definition}

These conditions are precisely the definition of a stable pseudo-flow (see Definition~\ref{def1}). This equivalence immediately leads to the following lemma:

\begin{lemma}
\label{lemma4}
The optimal solution of Programming~(\ref{eq5}) is a stable pseudo-flow.
\end{lemma}

Obviously, a stable pseudo-flow also satisfies the KKT conditions. By Lemma~\ref{lemma3}, we have
\begin{lemma}
\label{lemma5}
 A stable pseudo-flow  is the optimal solution of Programming~(\ref{eq5}).
\end{lemma}

By Lemma~\ref{lemma0}, Lemma~\ref{lemma4} and Lemma~\ref{lemma5}, Theorem~\ref{theo2} holds. 

\begin{theorem}
  \label{theo2}
  Given \(\{(s_k,t_k,d_k): k \in \mathcal{K}\}\) and edge capacities \(\{u_{ij}: (i,j)\in E\}\), the feasible region of Expression~(\ref{eq1}) is  empty if and only if there exists nonzero-stable pseudo-flow.
\end{theorem}

By Theorem~\ref{theo1} and Theorem~\ref{theo2}, the \emph{optimality condition} for  multicommodity flow problem may be summarized as follows:
\begin{enumerate}[label=(\roman*)]
\item if there exists a nonzero-stable pseudo-flow, there exists no feasible solution for the multicommodity flow problem; 
\item if there exists a zero-stable pseudo-flow, there exists a feasible solution and the zero-stable pseudo-flow is the feasible solution.
\end{enumerate}

\begin{remark}
In fact, both Theorem~\ref{theo1} and Theorem~\ref{theo2}  are true if the height function and 
the congestion function  satisfy the following conditions:

\begin{enumerate}[label=(\roman*)]
\item the height function \(h_{ik}(\Delta_{ik})\) is a strictly monotone increasing function and \(h_{ik}(0) = 0\); 
\item the congestion function \(\psi_{ij}\)   satisfies the following definition:
\begin{equation*}
\psi(f_{ij})=\left\{
\begin{aligned}
&0 ~~~~~~~~~~~~~~& if~f_{ij} \leq u_{ij} \\
&g(f_{ij} - u_{ij})~~~~~ & if~f_{ij} > u_{ij},
\end{aligned}
\right.
\end{equation*}
where \(g(0) = 0\) and \(g(\cdot)\) is a strictly monotone increasing function.
\end{enumerate}
\end{remark}

\begin{remark}
Expression~(\ref{cpgrad}) shows that the potential difference is exactly the same as the negative gradient of the objective function of Programming~(\ref{eq5}).
\end{remark}

\subsection{ Separable Structure}
\label{sec:alg}

The objective function of Programming (\ref{eq5}) has a mixed structure: the first term involves a sum over edges \((i,j)\), while the second term involves a sum over vertices \(i\). The flow variables \(f_{ij,k}\) are coupled in a complex way across both edges and vertices. We now rewrite its objective function as a single, unified sum over edges \((i,j)\).

Firstly, 
\begin{equation}\label{eq90}
\begin{aligned}
  &\int_{0}^{\Delta_{ik}}h_{ik}(\omega)d\omega = \frac{1}{2\delta_i}\Delta_{ik}^2 \\
&= \frac{1}{2\delta_i}( \sum_{j\in \delta^-(i)}f_{ji,k} - \sum_{j\in \delta^+(i)}f_{ij,k} + \hat{\Delta}_{ik})^2.\\
\end{aligned}
\end{equation}

By the definition of congestion function \(\psi_{ij}\), we have 
\begin{equation}\label{eq7}
\int_0^{f_{ij}}\psi_{ij}(\omega)d\omega =\left\{
\begin{aligned}
&0 ~~~~~~~~~~~~~~& if~f_{ij} \leq u_{ij} \\
&\frac{1}{2}(f_{ij} - u_{ij})^2~~~~~ & if~f_{ij} > u_{ij}.
\end{aligned}
\right.
\end{equation}

By introducing an auxiliary variable \(\{r_{ij}: r_{ij} \geq 0 \}\) which is the unused capacity for edge \((i,j)\), Programming~(\ref{eq5}) could be rewritten as 
\begin{equation}\label{eq8}
\begin{aligned}
&{\bf min}\quad z = \frac{1}{2}\sum_{(i,j) \in E} (\sum_{k\in \mathcal{K}} f_{ij,k} + r_{ij} - u_{ij})^2 \\
&+ \sum_{i,k}( \frac{1}{2\delta_i}( \sum_{j\in \delta^-(i)}f_{ji,k} - \sum_{j\in \delta^+(i)}f_{ij,k} + \hat{\Delta}_{ik})^2)\\
&= \frac{1}{2}\sum_{(i,j) \in E} [(\sum_{k\in \mathcal{K}} f_{ij,k} + r_{ij} - u_{ij})^2 \\
&+ \sum_{k\in \mathcal{K}}(\frac{\sum_{q\in \delta^-(i)}f_{qi,k} - \sum_{q\in \delta^+(i)}f_{iq,k} + \hat{\Delta}_{ik}}{\delta_{i}} )^2 \\
&+ \sum_{k\in \mathcal{K}}(\frac{\sum_{q\in \delta^-(j)}f_{qj,k} - \sum_{q\in \delta^+(j)}f_{jq,k} + \hat{\Delta}_{jk}}{\delta_{j}} )^2] \\
&= \frac{1}{2}\sum_{(i,j) \in E} [\psi_{ij}^2 + \sum_{k\in \mathcal{K}}h_{ik}^2 + \sum_{k\in \mathcal{K}}h_{jk}^2] \\
&{\bf s.t}\quad  f_{ij,k}\geq 0, \forall k\in \mathcal{K}, (i,j) \in E\\
\quad & r_{ij} \geq 0,\forall  (i,j) \in E \\
\end{aligned}
\end{equation}

The equivalence between the vertex-based and edge-based summations is established by a simple combinatorial principle. When a term is summed over all edges, any property related to the endpoints of those edges (like the height function \(h_{ik}\)) is counted for each vertex precisely as many times as its degree. In the expression \(\sum_{(i,j) \in E} [h_{ik}^2 + h_{jk}^2]\), a specific vertex \(v\) appears in the sum once for each edge connected to it. Therefore, the term \(h_{vk}^2\) is counted \(\delta_v\) times, justifying the transformation.

The objective function's edge-based terms are ideal for \emph{decomposition algorithms}, which split the problem into smaller, independent \emph{subproblems} per edge. These can be solved \emph{in parallel}, \emph{significantly} accelerating solutions for large networks.

\section{PDR via Exact Subproblem Optimization}

Theorem~\ref{theo1} and Theorem~\ref{theo2} above suggest a simple algorithmic approach for solving the multicommodity flow problem. That is, design an algorithm to obtain the stable pseudo-flow, which may be achieved  by adjusting the commodity flow  \(f_{ij,k}\) until the potential difference \(\phi_{ij,k} = 0\) or  \(f_{ij,k} = 0\). In this paper we call the algorithm the \emph{potential difference reduction} (PDR) algorithms.

The preceding section established that the multicommodity flow problem can be formulated as the convex optimization problem in Expression~(\ref{eq8}), 
which seeks a zero-stable pseudo-flow. The  parallel algorithm proposed here is to iteratively adjust the flows on each edge to locally minimize the objective function.

\subsection{The Edge Subproblem}

For each edge \((i,j) \in E \), we define an edge subproblem, which is derived from the terms in the global objective function (Expression~\ref{eq8}).
The objective is to minimize a function representing the new congestion on edge \((i,j)\) plus the change in "height energy" at vertices \(i\) and \(j\) 
caused by the flow update on that single edge. This yields the following Quadratic Program (QP):

\begin{equation}\label{eq9}
  \begin{aligned}
  {\bf min} \quad &z_{ij} = \frac{1}{2}[(\sum_{k\in \mathcal{K}}f_{ij,k} +  \sum_{k\in \mathcal{K}}\Delta f_{ij,k} + r_{ij} - u_{ij})^2 \\
  &+ \sum_{k\in \mathcal{K}}(h_{ik}-\Delta f_{ij,k})^2 + \sum_{k\in \mathcal{K}}(h_{jk}+\Delta f_{ij,k})^2]\\
  {\bf s.t}\quad &  \Delta f_{ij,k} \geq -f_{ij,k}, \forall k\in \mathcal{K}\\
  \quad & r_{ij} \geq 0 \\
  \end{aligned}
  \end{equation}

These terms, $(h_{ik}-\Delta f_{ij,k})$ and $(h_{jk}+\Delta f_{ij,k})$, reflect the local impact of the flow adjustment on node heights: as a flow increment $\Delta f_{ij,k}$ moves from node $i$ to node $j$, it effectively lowers the height at $i$ while raising the height at $j$.

Each edge subproblem is solved in parallel, based on the height information from the previous iteration. Note that it does not account for the simultaneous flow changes happening on other edges connected to vertices \(i\) and \(j\). 

\begin{remark}
Since this subproblem is defined for a \textbf{single edge}, its local contribution to the height update is treated as if the degree were 1. This explains why the flow change $\Delta f_{ij,k}$ is subtracted directly from the height $h_{ik}$, rather than being scaled by the total vertex degree $\delta_i$.
\end{remark}
  
\subsection{Algorithm Description}

The proposed algorithm, which we term \textbf{PDR via Exact Subproblem Optimization}, iteratively solves these edge subproblems  until the system converges to a stable pseudo-flow. The steps are outlined in Algorithm~\ref{alg:qp_projection}.

Each iteration begins with parallel computations of the flow imbalance \(\Delta_{ik}\), vertex heights \(h_{ik}\), and edge congestions \(\psi_{ij}\). The algorithm then checks for convergence by determining if the change in the objective function value between iterations is within a specified tolerance \(\epsilon_{obj}\). If not, it proceeds to the main update step, where it solves the QP edge subproblem in Expression~(\ref{eq10}) for each edge \((i, j)\) to find the flows for the next iteration, \(f_{ij,k}^{(t+1)}\). Note that Expression~(\ref{eq10}) is simply a variable substitution of Expression~(\ref{eq9}).

Crucially, since each subproblem depends only on local information, these optimizations can be executed \textbf{in parallel} for all edges. This parallel structure makes the algorithm exceptionally scalable and well-suited for large networks.

\begin{equation}\label{eq10}
    \begin{aligned}
    {\bf min}\quad &z_{ij} = \frac{1}{2} [(\sum_{k\in \mathcal{K}} f_{ij,k}^{(t+1)} + r_{ij} - u_{ij})^2 \\
    &+ \sum_{k\in \mathcal{K}}(h_{ik}^{(t)} + f_{ij,k}^{(t)} - f_{ij,k}^{(t+1)})^2 \\
    &+ \sum_{k\in \mathcal{K}}(h_{jk}^{(t)} - f_{ij,k}^{(t)} + f_{ij,k}^{(t+1)})^2]\\
    {\bf s.t}\quad & f_{ij,k}^{(t+1)}\geq 0, \forall k\in \mathcal{K}\\
    \quad & r_{ij} \geq 0 \\
    \end{aligned}
    \end{equation}

    \begin{algorithm}[H]
        \caption{PDR via Exact Subproblem Optimization for MCF}
        \label{alg:qp_projection}
        \begin{algorithmic}[1]
        \State \textbf{Input:} Graph $G(V,E)$, capacities $\{u_{ij}\}$, commodities $\{(s_k, t_k, d_k)\}$, tolerance $\epsilon_{\text{obj}} > 0$.
        \State \textbf{Initialize:} Set initial flow $f^{(0)} \leftarrow 0$, iteration $t \leftarrow 0$, previous objective value $z_{\text{prev}} \leftarrow \infty$.
        
        \While{$t < \text{max\_iterations}$}
            \ForAll{vertices $i \in V$, commodities $k \in \mathcal{K}$ \textbf{in parallel}}
                \State Flow imbalance $\Delta_{ik}^{(t)} \leftarrow \sum_{j \in \delta^-(i)} f_{ji,k}^{(t)} - \sum_{j \in \delta^+(i)} f_{ij,k}^{(t)} + \hat{\Delta}_{ik}$
                \State Vertex height $h_{ik}^{(t)} \leftarrow \Delta_{ik}^{(t)} / \delta_i$
            \EndFor
            \ForAll{edges $(i,j) \in E$ \textbf{in parallel}}
                \State Edge congestion $\psi_{ij}^{(t)} \leftarrow \max\left(0, \sum_{k \in \mathcal{K}} f_{ij,k}^{(t)} - u_{ij}\right)$
            \EndFor
            
            \State \Comment{Compute current objective value $z^{(t)}$}
            \State $z^{(t)} \leftarrow \frac{1}{2} \sum_{(i,j) \in E} (\psi_{ij}^{(t)})^2 + \frac{1}{2} \sum_{i \in V, k \in \mathcal{K}} \delta_i (h_{ik}^{(t)})^2$
            
            \State \Comment{Check for convergence based on objective value change}
            \If{$|z^{(t)} - z_{\text{prev}}| < \epsilon_{\text{obj}}$}
                \State \textbf{break}
            \EndIf
            
            \State $z_{\text{prev}} \leftarrow z^{(t)}$ \Comment{Update previous objective value for the next iteration}
            
            \State \Comment{Update flows by solving edge subproblems}
            \State For each edge $(i,j) \in E$ \textbf{in parallel}, solve the QP subproblem in Expression~(14) to obtain the updated flows $f_{ij,k}^{(t+1)}$ for all $k \in \mathcal{K}$.
            
            \State $t \leftarrow t + 1$
        \EndWhile
        \State \textbf{Return} final flow $f^{(t)}$
        \end{algorithmic}
    \end{algorithm}

	\subsection{Proof of Algorithm Convergence}

	In this section, we prove that Algorithm~\ref{alg:qp_projection}, is convergent. 
	The core idea is to show that the value of the global objective function $z$, defined in Expression~(\ref{eq8}), is monotonically non-increasing in each iteration. 
	Since $z$ is a sum of squared terms, it is clearly bounded below ($z \ge 0$). 
	A monotonically non-increasing sequence that is bounded below is guaranteed to converge.

	\begin{theorem}
	For the sequence of flows $\{ \boldsymbol{f}^{(t)} \}_{t=0,1,2,\cdots}$ generated by Algorithm~\ref{alg:qp_projection}, the corresponding sequence of global objective function values $\{ z(\boldsymbol{f}^{(t)}) \}$ is monotonically non-increasing. That is:
	$$z(\boldsymbol{f}^{(t+1)}) \le z(\boldsymbol{f}^{(t)}), \quad \forall t \ge 0$$
	\end{theorem}
	
	\proof
	To prove this theorem, we introduce a surrogate function and complete the proof in two steps.
	
	Let $\boldsymbol{f}^{(t)}$ be the flow vector at the start of iteration $t$. The algorithm computes the new flow vector $\boldsymbol{f}^{(t+1)}$ by solving, in parallel, the QP subproblems (Expression~(\ref{eq10})) for all edges $(i,j) \in E$. The sum of the objective functions of all these subproblems can be formulated as a surrogate function $Q(\boldsymbol{f}', \boldsymbol{f}^{(t)})$, where $\boldsymbol{f}'$ is the variable to be optimized. Based on Expression~(\ref{eq10}), we define $Q$ as:
	
	$$Q(\boldsymbol{f}', \boldsymbol{f}^{(t)}) = \frac{1}{2}\sum_{(i,j) \in E} \left[ \left(\sum_{k\in \mathcal{K}} f'_{ij,k} + r'_{ij} - u_{ij}\right)^2 + \sum_{k\in \mathcal{K}}\left(h_{ik}^{(t)} - (f'_{ij,k} - f_{ij,k}^{(t)})\right)^2 + \sum_{k\in \mathcal{K}}\left(h_{jk}^{(t)} + (f'_{ij,k} - f_{ij,k}^{(t)})\right)^2 \right]$$
	
	Here, $h_{ik}^{(t)}$ and $h_{jk}^{(t)}$ are the vertex heights calculated based on the flow $\boldsymbol{f}^{(t)}$.
	
	\emph{Step 1: Proving that $Q(\boldsymbol{f}^{(t+1)}, \boldsymbol{f}^{(t)}) \le z(\boldsymbol{f}^{(t)})$}
	
	By the definition of Algorithm~\ref{alg:qp_projection}, $\boldsymbol{f}^{(t+1)}$ is obtained by minimizing $Q(\boldsymbol{f}', \boldsymbol{f}^{(t)})$, as this minimization is composed of the independent minimizations of each edge subproblem. Therefore, $\boldsymbol{f}^{(t+1)}$ is the global minimizer of $Q(\boldsymbol{f}', \boldsymbol{f}^{(t)})$.
	$$\boldsymbol{f}^{(t+1)} = \arg\min_{\boldsymbol{f}' \ge 0} Q(\boldsymbol{f}', \boldsymbol{f}^{(t)})$$
	This implies that for any other feasible flow, including $\boldsymbol{f}^{(t)}$ itself, we have $Q(\boldsymbol{f}^{(t+1)}, \boldsymbol{f}^{(t)}) \le Q(\boldsymbol{f}^{(t)}, \boldsymbol{f}^{(t)})$.
	
	Now, let's evaluate $Q(\boldsymbol{f}^{(t)}, \boldsymbol{f}^{(t)})$. When $\boldsymbol{f}' = \boldsymbol{f}^{(t)}$, the change in flow is zero ($f'_{ij,k} - f_{ij,k}^{(t)} = 0$). Substituting this into the definition of $Q$ yields:
	$$Q(\boldsymbol{f}^{(t)}, \boldsymbol{f}^{(t)}) = \frac{1}{2}\sum_{(i,j) \in E} \left[ \left(\sum_{k\in \mathcal{K}} f_{ij,k}^{(t)} + r_{ij}^{(t)} - u_{ij}\right)^2 + \sum_{k\in \mathcal{K}}(h_{ik}^{(t)})^2 + \sum_{k\in \mathcal{K}}(h_{jk}^{(t)})^2 \right]$$
	This is identical to the definition of the global objective function $z(\boldsymbol{f}^{(t)})$ in Expression~(\ref{eq8}). Thus, we have $Q(\boldsymbol{f}^{(t)}, \boldsymbol{f}^{(t)}) = z(\boldsymbol{f}^{(t)})$.
	Combining these facts gives us our first key inequality:
	$$Q(\boldsymbol{f}^{(t+1)}, \boldsymbol{f}^{(t)}) \le z(\boldsymbol{f}^{(t)})$$
	
	\emph{Step 2: Proving that $z(\boldsymbol{f}^{(t+1)}) \le Q(\boldsymbol{f}^{(t+1)}, \boldsymbol{f}^{(t)})$}
	
	Next, we compare the true global objective value $z(\boldsymbol{f}^{(t+1)})$ under the new flow with the value of the surrogate function $Q(\boldsymbol{f}^{(t+1)}, \boldsymbol{f}^{(t)})$.
	$$z(\boldsymbol{f}^{(t+1)}) = \frac{1}{2}\sum_{(i,j) \in E} \left[ (\psi_{ij}^{(t+1)})^2 + \sum_{k\in \mathcal{K}}(h_{ik}^{(t+1)})^2 + \sum_{k\in \mathcal{K}}(h_{jk}^{(t+1)})^2 \right]$$
	$$Q(\boldsymbol{f}^{(t+1)}, \boldsymbol{f}^{(t)}) = \frac{1}{2}\sum_{(i,j) \in E} \left[ \psi_{ij}^{(t+1)2} + \sum_{k\in \mathcal{K}}(h_{ik}^{(t)} - \Delta f_{ij,k})^2 + \sum_{k\in \mathcal{K}}(h_{jk}^{(t)} + \Delta f_{ij,k})^2 \right]$$
	where $\Delta f_{ij,k} = f_{ij,k}^{(t+1)} - f_{ij,k}^{(t)}$.
	
	Comparing the two expressions, the congestion terms $(\psi_{ij}^{(t+1)})^2$ are identical. We only need to compare the height terms. 
	The height at vertex $i$ for commodity $k$ under the new flow $\boldsymbol{f}^{(t+1)}$ is given by:
		
	$$h_{ik}^{(t+1)} = h_{ik}^{(t)} + \frac{1}{\delta_i}\left( \sum_{j\in \delta^-(i)}\Delta f_{ji,k} - \sum_{j\in \delta^+(i)}\Delta f_{ij,k} \right).$$
	
	To apply Jensen's inequality, we rewrite $h_{ik}^{(t+1)}$ as an average by distributing $h_{ik}^{(t)} = \frac{1}{\delta_i}\sum_{j \in \delta(i)} h_{ik}^{(t)}$:
	
	$$h_{ik}^{(t+1)}  = \frac{1}{\delta_i} \left( \sum_{j \in \delta^+(i)} (h_{ik}^{(t)} - \Delta f_{ij,k}) + \sum_{j \in \delta^-(i)} (h_{ik}^{(t)} + \Delta f_{ji,k}) \right)$$
	
	Since the function $x^2$ is convex, we can apply Jensen's inequality, which states that the square of an average is less than or equal to the average of the squares:
	
	$$(h_{ik}^{(t+1)})^2 \le \frac{1}{\delta_i} \left( \sum_{j \in \delta^+(i)} (h_{ik}^{(t)} - \Delta f_{ij,k})^2 + \sum_{j \in \delta^-(i)} (h_{ik}^{(t)} + \Delta f_{ji,k})^2 \right)$$
	
	Multiplying by $\delta_i$ and summing over all vertices $i$ and commodities $k$:
	
	$$\sum_{i \in V, k \in K} \delta_i (h_{ik}^{(t+1)})^2 \le \sum_{i \in V, k \in K} \left( \sum_{j \in \delta^+(i)} (h_{ik}^{(t)} - \Delta f_{ij,k})^2 + \sum_{j \in \delta^-(i)} (h_{ik}^{(t)} + \Delta f_{ji,k})^2 \right)$$
	
	The left-hand side (LHS) is the height-energy component of $2 \cdot z(\boldsymbol{f}^{(t+1)})$. We can rearrange the right-hand side (RHS) by summing over edges $(i,j)$ instead of vertices. Each edge $(i,j)$ contributes exactly two terms to the sum: one for vertex $i$ (as an outgoing edge) and one for vertex $j$ (as an incoming edge). This gives:
	
	$$\text{RHS} = \sum_{(i,j) \in E} \sum_{k \in K} \left[ (h_{ik}^{(t)} - \Delta f_{ij,k})^2 + (h_{jk}^{(t)} + \Delta f_{ij,k})^2 \right]$$
	
	This is precisely the height-energy component of $2 \cdot Q(\boldsymbol{f}^{(t+1)}, \boldsymbol{f}^{(t)})$. Since the congestion terms in $z(\boldsymbol{f}^{(t+1)})$ and $Q(\boldsymbol{f}^{(t+1)}, \boldsymbol{f}^{(t)})$ are identical, this establishes our second key inequality:
	
	$$z(\boldsymbol{f}^{(t+1)}) \le Q(\boldsymbol{f}^{(t+1)}, \boldsymbol{f}^{(t)})$$

	Combining our two inequalities, we have:
	$$z(\boldsymbol{f}^{(t+1)}) \le Q(\boldsymbol{f}^{(t+1)}, \boldsymbol{f}^{(t)}) \le z(\boldsymbol{f}^{(t)})$$
	This proves that the sequence of global objective function values $\{ z(\boldsymbol{f}^{(t)}) \}$ is monotonically non-increasing. Since $z(\boldsymbol{f}) \ge 0$, the sequence is bounded below. By the Monotone Convergence Theorem, the sequence must converge to a limit. This implies that the flow vector $\boldsymbol{f}^{(t)}$ converges to a fixed point, which, by its construction as a solution to the KKT conditions, is a stable pseudo-flow.

	\endproof

	\subsection{Exact Solution of the Edge Subproblem}
	
	The effectiveness of the PDR algorithm hinges on our ability to solve the edge subproblem (Expression~\ref{eq10}) efficiently and exactly. While it is a Quadratic Program (QP), its specific structure allows for a closed-form analytical solution. This section derives this solution by applying the Karush-Kuhn-Tucker (KKT) conditions, demonstrating that the optimization for each edge reduces to a simple, fast computation.
	
	First, let's restate the objective function for a single edge $(i,j)$ at iteration $t$. We seek to find the next-iteration flows, which we'll denote as $f'_{ij,k}$ for clarity, that minimize $z_{ij}$:
	\begin{equation*}
		\begin{aligned}
		\min_{f'_{ij,k} \ge 0, r'_{ij} \ge 0} \quad &z_{ij} = \frac{1}{2} \left(\sum_{k \in \mathcal{K}} f'_{ij,k} + r'_{ij} - u_{ij}\right)^2 \\
		&+ \frac{1}{2}\sum_{k \in \mathcal{K}}\left(h_{ik}^{(t)} + f_{ij,k}^{(t)} - f'_{ij,k}\right)^2 \\
		&+ \frac{1}{2}\sum_{k \in \mathcal{K}}\left(h_{jk}^{(t)} - f_{ij,k}^{(t)} + f'_{ij,k}\right)^2 \\
		\end{aligned}
	\end{equation*}
	This objective function is strictly convex and its constraints are linear, so the KKT conditions are both necessary and sufficient for optimality.

	\subsubsection{Karush-Kuhn-Tucker (KKT) Formulation}
	
	We introduce Lagrange multipliers $\lambda_k \ge 0$ for the non-negativity constraints $f'_{ij,k} \ge 0$ and $\eta \ge 0$ for $r'_{ij} \ge 0$. The Lagrangian $\mathcal{L}$ is:
	$$
	\mathcal{L} = z_{ij} - \sum_{k \in \mathcal{K}} \lambda_k f'_{ij,k} - \eta r'_{ij}
	$$
	The stationarity condition requires the gradient of $\mathcal{L}$ with respect to the primal variables ($f'_{ij,k}$ and $r'_{ij}$) to be zero.
	
	\begin{enumerate}
		\item \textbf{Derivative with respect to $r'_{ij}$:} Let $\psi'_{ij} = \sum_{k \in \mathcal{K}} f'_{ij,k} + r'_{ij} - u_{ij}$ be the new congestion value.
		$$
		\frac{\partial \mathcal{L}}{\partial r'_{ij}} = \left(\sum_{k \in \mathcal{K}} f'_{ij,k} + r'_{ij} - u_{ij}\right) - \eta = \psi'_{ij} - \eta = 0 \implies \boldsymbol{\psi'_{ij} = \eta}
		$$
		This provides a direct physical interpretation: the optimal congestion level $\psi'_{ij}$ is equal to the dual variable $\eta$ associated with the slack capacity.
	
		\item \textbf{Derivative with respect to $f'_{ij,k}$:}
		$$
		\frac{\partial \mathcal{L}}{\partial f'_{ij,k}} = \psi'_{ij} - (h_{ik}^{(t)} + f_{ij,k}^{(t)} - f'_{ij,k}) + (h_{jk}^{(t)} - f_{ij,k}^{(t)} + f'_{ij,k}) - \lambda_k = 0
		$$
		Rearranging the terms to solve for $f'_{ij,k}$:
		$$
		\psi'_{ij} + 2f'_{ij,k} - (h_{ik}^{(t)} - h_{jk}^{(t)} + 2f_{ij,k}^{(t)}) - \lambda_k = 0
		$$
		To simplify, we define a constant "driving potential" $c_k$ for each commodity $k$, which combines all known values from iteration $t$:
		$$ c_k \equiv \frac{1}{2}(h_{ik}^{(t)} - h_{jk}^{(t)} + 2f_{ij,k}^{(t)} )$$
		The stationarity condition then simplifies to $2f'_{ij,k} = 2c_k - \psi'_{ij} + \lambda_k$.
	\end{enumerate}
	
	\subsubsection{Deriving the Flow Update Rule}
	
	The complementary slackness condition, $\lambda_k f'_{ij,k} = 0$, allows us to determine the optimal flow $f'_{ij,k}$:
	\begin{itemize}
		\item If $f'_{ij,k} > 0$: Complementary slackness requires $\lambda_k = 0$. The flow is $f'_{ij,k} = c_k - \frac{\psi'_{ij}}{2}$. This is only valid if $2c_k > \psi'_{ij}$.
		\item If $f'_{ij,k} = 0$: $\lambda_k \ge 0$ is possible. The stationarity condition gives $\lambda_k = \psi'_{ij} - 2c_k \ge 0$, which means this case holds if $2c_k \le \psi'_{ij}$.
	\end{itemize}
	Combining these two outcomes gives a single, elegant expression for the optimal flow as a function of the yet-unknown congestion level $\psi'_{ij}$:
	\begin{equation} \label{eq:flow_update_psi}
	f'_{ij,k} = \max\left(0, c_k - \frac{\psi'_{ij}}{2}\right)
	\end{equation}
	The problem has been reduced to finding a single scalar value, $\psi'_{ij}$.
	
	\subsubsection{Solving for the Congestion Level $\psi'_{ij}$}
	
	The final KKT condition, $\eta r'_{ij} = 0$, combined with our finding that $\eta = \psi'_{ij}$, gives $\psi'_{ij} r'_{ij} = 0$. This implies that either the congestion is zero or the slack capacity is zero, leading to two distinct cases.
	
	\paragraph{Case 1: Uncongested Edge ($r'_{ij} > 0$)}
	If the optimal solution has unused capacity, then we must have $\boldsymbol{\psi'_{ij} = 0}$. Substituting this into Equation~(\ref{eq:flow_update_psi}) gives the candidate solution:
	$$
	f'_{ij,k} = \max\left(0, c_k\right)
	$$
	This solution is valid if and only if the resulting total flow is less than the capacity, which is the condition for $r'_{ij} > 0$ to be possible. We check this condition:
	$$
	\text{If} \quad \sum_{k \in \mathcal{K}} \max\left(0, c_k\right) < u_{ij}, \quad \text{then this is the optimal solution.}
	$$

	\paragraph{Case 2: Congested Edge ($r'_{ij} = 0$)}
	If the check for Case 1 fails (i.e., $\sum_k \max(0, c_k) \ge u_{ij}$), the edge must be at or over capacity, meaning the optimal slack is $\boldsymbol{r'_{ij} = 0}$. In this case, the congestion $\psi'_{ij}$ can be positive.
	
	With $r'_{ij} = 0$, the definition of $\psi'_{ij}$ becomes $\psi'_{ij} = \sum_k f'_{ij,k} - u_{ij}$. We find the value of $\psi'_{ij}$ by substituting our expression for $f'_{ij,k}$ (Equation~\ref{eq:flow_update_psi}) into this definition:
	$$
	\psi'_{ij} = \sum_{k \in \mathcal{K}} \max\left(0, c_k - \frac{\psi'_{ij}}{2}\right) - u_{ij}
	$$
	The right-hand side of this equation is a piecewise-linear, continuous, and monotonically decreasing function of $\psi'_{ij}$. We are looking for the unique root of this equation. This root can be found very efficiently by first sorting the potential values $c_k$. The breakpoints of the function occur at $\psi'_{ij} = 2c_k$. By iterating through the sorted list of potentials, we can quickly identify the linear segment that contains the root. This procedure takes $O(K \log K)$ time, dominated by sorting the $K$ commodity potentials. 
	
	Once the unique root $\psi'_{ij} \ge 0$ is found, the final optimal flows are calculated using Equation~(\ref{eq:flow_update_psi}). The complete algorithm for this procedure is detailed in Algorithm~\ref{alg:qp-solver-final}.

	\begin{algorithm}[H]
		\caption{Solver for the Edge Subproblem}
		\label{alg:qp-solver-final}
		\begin{algorithmic}[1]
		\Require Current vertex heights and flows $\{h_{ik}^t, h_{jk}^t, f_{ij,k}^{t}\}_{k \in \mathcal{K}}$, and edge capacity $u_{ij}$.
		\Ensure Optimal next-iteration flows $f_{ij,k}^{t+1}$, slack capacity $r_{ij}$, and new congestion $\psi_{ij}$.
		
		\State \textbf{1. Preprocessing:}
		\State For each commodity $k \in \mathcal{K}$, compute its potential:
		\Statex \qquad $c_k \gets \frac{1}{2}\left( (h_{ik}^t + f_{ij,k}^{t}) - (h_{jk}^t - f_{ij,k}^{t}) \right)$.
		
		\State \textbf{2. Check for Congestion:}
		\State Compute the potential flow sum $s \gets \sum_{k \in \mathcal{K}} \max(0, c_k)$.
		
		\If{$s < u_{ij}$} \Comment{\textit{Case 1: The edge is not congested.}}
			\State $r_{ij} \gets u_{ij} - s$.
			\State $\psi_{ij} \gets 0$.
		\Else \Comment{\textit{Case 2: The edge is congested, solve for $\psi_{ij}$.}}
			\State $r_{ij} \gets 0$.
			\State Let $\mathcal{K}^+ = \{k \in \mathcal{K} \mid c_k > 0\}$.
			\State Sort potentials in $\mathcal{K}^+$ in descending order: $c_{(1)} \ge c_{(2)} \ge \dots \ge c_{(m)}$.
			\State $\text{sum\_c} \gets 0$.
			\For{$p = 1, \dots, m$}
				\State $\text{sum\_c} \gets \text{sum\_c} + c_{(p)}$.
				\State $\psi_{\text{candidate}} \gets (\text{sum\_c} - u_{ij}) / (1 + p/2)$.
		
				\State \Comment{A candidate solution is valid only if it falls within the assumed interval.}
				\State Let $c_{\text{upper\_bound}} = c_{(p)}$.
				\State Let $c_{\text{lower\_bound}} = (p < m) ? c_{(p+1)} : 0$.
		
				\If{$2 \cdot c_{\text{lower\_bound}} < \psi_{\text{candidate}} \le 2 \cdot c_{\text{upper\_bound}}$}
					\State $\psi_{ij} \gets \psi_{\text{candidate}}$. \Comment{This is the unique, correct root.}
					\State \textbf{break} \Comment{Exit the loop.}
				\EndIf
			\EndFor
		\EndIf
		
		\State \textbf{3. Final Flow Calculation:}
		\State For each commodity $k \in \mathcal{K}$:
		\Statex \qquad $f_{ij,k}^{t+1} \gets \max(0, c_k - \psi_{ij}/2)$.
		
		\State \Return $(f_{ij,k}^{t+1})_{k \in \mathcal{K}}$, $r_{ij}$, and $\psi_{ij}$.
		\end{algorithmic}
		\end{algorithm}

	\section{PDR via Heuristic Algorithm}
	
	While the PDR algorithm with an exact subproblem solver guarantees convergence, its practical performance can be limited when the number of commodities, $|\mathcal{K}|$, is very large. The efficiency of the exact solver (Algorithm \ref{alg:qp-solver-final}) depends on sorting the commodity potentials, an operation with a computational complexity of $O(|\mathcal{K}| \log |\mathcal{K}|)$ for each edge. This sorting step can become a significant computational bottleneck in large-scale scenarios.
	To address this challenge, we propose more computationally efficient heuristic variants based on the \textbf{gradient projection} method.
	
	\subsection{Adaptive Gradient Descent Algorithm}
	
	This approach leverages a key insight from our formulation: as noted in Remark 2, the negative gradient of the global objective function $z$ with respect to the flow $f_{ij,k}$ is precisely the potential difference, $\phi_{ij,k}$. The potential difference therefore indicates the direction of steepest descent for the objective function. This naturally leads to a gradient descent heuristic, which adjusts the flow on each edge in proportion to this gradient.
	
	The core of the algorithm lies in the flow update rule. For each edge $(i, j)$ and commodity $k$, a proposed flow change is computed based on the potential difference and a per-edge \textbf{learning rate} $\beta_{ij}$. This change is then \textbf{projected} onto the feasible set to ensure $f_{ij,k}^{(t+1)} \geq 0$:
	$$ \phi_{ij,k}^{(t)} = (h_{ik}^{(t)} - h_{jk}^{(t)} - \psi_{ij}^{(t)})$$
	$$\Delta f_{ij,k} \leftarrow \max(-f_{ij,k}^{(t)},\beta_{ij} \cdot \phi_{ij,k}^{(t)})$$

	Before applying an update, the algorithm first evaluates whether the step $\{\Delta f_{ij,k}\}_{k \in \mathcal{K}}$ for an edge $(i,j)$ will decrease the \textbf{local objective function} $z_{ij}$ of Expression~\ref{eq9}. If the update would lead to an increase in the objective value, the step for that edge is rejected for the current iteration. This acceptance/rejection criterion ensures that the objective function is monotonically non-increasing, which is a key property for ensuring convergence.
	
	The learning rate $\beta_{ij}$ for each edge is adapted dynamically. If a step is rejected, it indicates the learning rate was too aggressive, so it is multiplicatively decreased. Conversely, the algorithm periodically attempts to increase the learning rate to accelerate progress. This adaptive scheme allows the algorithm to be both aggressive in its search for the minimum and cautious when necessary, improving overall performance.

	The complete procedure, outlined in Algorithm~\ref{alg:gradient_descent}, provides a computationally inexpensive alternative to the exact subproblem solver. It trades guaranteed per-iteration optimality for extremely fast, scalable iterations, making it particularly well-suited for problems with a very large number of commodities.

	\begin{algorithm}[H]
		\caption{PDR via Adaptive Gradient Descent}
		\label{alg:gradient_descent}
		\small
		\setlength{\itemsep}{0pt} 
		\setlength{\parsep}{0pt}    
		\begin{algorithmic}[1]
		\State \textbf{Input:} Graph $G(V,E)$, capacities $\{u_{ij}\}$, commodities $\{(s_k, t_k, d_k)\}$, tolerance $\epsilon_{\text{obj}} > 0$.
		\State \textbf{Parameters:} Initial learning rate $\beta_{init}$, min/max rates $\beta_{min}, \beta_{max}$, factors $\alpha_{inc} > 1$, $\alpha_{dec} < 1$.
		\State \textbf{Initialize:} Set initial flow  $f^{(0)} \leftarrow 0$, iteration $t \leftarrow 0$, learning rate $\beta_{ij} \leftarrow \beta_{init}, \forall (i,j) \in E$, previous objective value $z_{\text{prev}} \leftarrow \infty$.
		
		\While{$t < \text{max\_iterations}$}
			\ForAll{vertices $i \in V$, commodities $k \in \mathcal{K}$ \textbf{in parallel}}
				\State Flow imbalance $\Delta_{ik}^{(t)} \leftarrow \sum_{j \in \delta^-(i)} f_{ji,k}^{(t)} - \sum_{j \in \delta^+(i)} f_{ij,k}^{(t)} + \hat{\Delta}_{ik}$
				\State Vertex height $h_{ik}^{(t)} \leftarrow \Delta_{ik}^{(t)} / \delta_i$
			\EndFor
			\ForAll{edges $(i,j) \in E$ \textbf{in parallel}}
				\State Edge congestion $\psi_{ij}^{(t)} \leftarrow \max\left(0, \sum_{k \in \mathcal{K}} f_{ij,k}^{(t)} - u_{ij}\right)$
			\EndFor
			\State $z^{(t)} \leftarrow \frac{1}{2} \sum_{(i,j) \in E} (\psi_{ij}^{(t)})^2 + \frac{1}{2} \sum_{i \in V, k \in \mathcal{K}} \delta_i (h_{ik}^{(t)})^2$
			\If {$|z^{(t)} - z_{\text{prev}}| < \epsilon_{\text{obj}}$}
			\State \textbf{break}
			\EndIf
			\State $z_{\text{prev}} \leftarrow z^{(t)}$ 
			\State \Comment{Update flows using a projected gradient step}
			\ForAll{edges $(i,j) \in E$ \textbf{in parallel}}
				\State For each commodity $k \in \mathcal{K}$, $\phi_{ij,k} \leftarrow h_{ik}^{(t)} - h_{jk}^{(t)} - \psi_{ij}^{(t)}$.
				\State For each commodity $k \in \mathcal{K}$, $\Delta f_{ij,k} \leftarrow \max(-f_{ij,k}^{(t)}, \beta_{ij} \cdot \phi_{ij,k})$.
				\State Estimate objective change $\Delta z_{ij}$ if update $\{\Delta f_{ij,k}\}$ is applied.
				
				\If{$\Delta z_{ij} \le 0$} 
					\State For each commodity $k \in \mathcal{K}$, $f_{ij,k}^{(t+1)} \leftarrow f_{ij,k}^{(t)} + \Delta f_{ij,k}$.
				\Else \Comment{Reject update and decrease learning rate}
					\State For each commodity $k \in \mathcal{K}$, $f_{ij,k}^{(t+1)} \leftarrow f_{ij,k}^{(t)}$.
					\State $\beta_{ij} \leftarrow \max(\beta_{min}, \beta_{ij} \cdot \alpha_{dec})$.
				\EndIf
				\If{$t \pmod{10} = 0$} 
					\State $\beta_e \leftarrow \beta_e \cdot \alpha_{inc}$ \Comment{Periodically increase learning rate}
				\EndIf
			\EndFor
			
			\State $t \leftarrow t + 1$
		\EndWhile
		\State \textbf{Return} final flow $f^{(t)}$
		\end{algorithmic}
	\end{algorithm}

	\subsection{Momentum-Based Adaptive Gradient Descent Algorithm}
	
	While the Adaptive Gradient Descent method provides a significant computational advantage, there are numerous techniques to accelerate gradient-based methods. Here, we introduce a simple yet powerful and widely-used approach: incorporating \textbf{momentum}. The resulting method, the Momentum-Based Adaptive Gradient Descent algorithm, is detailed in Algorithm~\ref{alg:momentum_gradient}.

	To implement this, we introduce a \textbf{velocity vector}, $\boldsymbol{v}$, for each commodity on each edge. This vector accumulates a history of the gradient-based steps, effectively smoothing the optimization trajectory. As shown in Algorithm~\ref{alg:momentum_gradient}, the update process is modified as follows.   First, the standard projected  step, $\Delta f_{ij,k}$, is calculated as in the previous method. Then, the velocity is updated by blending its previous value with the new step, controlled by a momentum coefficient, $\gamma \in [0, 1)$:
	$$v_{ij,k}^{(t+1)} \leftarrow \gamma \cdot v_{ij,k}^{(t)} + \Delta f_{ij,k}$$
	The flow is then updated by applying this new velocity, and the result is projected to ensure non-negativity: $f_{ij,k}^{(t+1)} \leftarrow \max(0, f_{ij,k}^{(t)} + v_{ij,k}^{(t+1)})$.

	Note that Algorithm \ref{alg:momentum_gradient} is identical to Algorithm \ref{alg:gradient_descent}, with the exception of the modifications highlighted in \textcolor{blue}{blue}. Numerical experiments show that such a small modification can significantly accelerate convergence.

	\begin{algorithm}[H]
		\caption{PDR via Gradient Descent with Momentum}
		\label{alg:momentum_gradient}
		\small
		\setlength{\itemsep}{0pt} 
		\setlength{\parsep}{0pt}   
		\begin{algorithmic}[1]
		\State \textbf{Input:} Graph $G(V,E)$, capacities $\{u_{ij}\}$, commodities $\{(s_k, t_k, d_k)\}$, tolerance $\epsilon_{\text{obj}} > 0$.
		\State \textbf{Parameters:} Initial learning rate $\beta_{init}$, min/max rates $\beta_{min}, \beta_{max}$,  factors $\alpha_{inc} > 1$, $\alpha_{dec} < 1$, momentum coefficient $\gamma \in [0, 1)$.
		\State \textbf{Initialize:} Set initial flow  $f^{(0)} \leftarrow 0$, velocity $v^{(0)} \leftarrow 0$, iteration $t \leftarrow 0$, learning rate $\beta_{ij} \leftarrow \beta_{init}, \forall (i,j) \in E$.
		
		\While{$t < \text{max\_iterations}$}
			\ForAll{vertices $i \in V$, commodities $k \in \mathcal{K}$ \textbf{in parallel}}
				\State Flow imbalance $\Delta_{ik}^{(t)} \leftarrow \sum_{j \in \delta^-(i)} f_{ji,k}^{(t)} - \sum_{j \in \delta^+(i)} f_{ij,k}^{(t)} + \hat{\Delta}_{ik}$
				\State Vertex height $h_{ik}^{(t)} \leftarrow \Delta_{ik}^{(t)} / \delta_i$
			\EndFor
			\ForAll{edges $(i,j) \in E$ \textbf{in parallel}}
				\State Edge congestion $\psi_{ij}^{(t)} \leftarrow \max\left(0, \sum_{k \in \mathcal{K}} f_{ij,k}^{(t)} - u_{ij}\right)$
			\EndFor
			\State $z^{(t)} \leftarrow \frac{1}{2} \sum_{(i,j) \in E} (\psi_{ij}^{(t)})^2 + \frac{1}{2} \sum_{i \in V, k \in \mathcal{K}} \delta_i (h_{ik}^{(t)})^2$
			\If {$|z^{(t)} - z_{\text{prev}}| < \epsilon_{\text{obj}}$}
			\State \textbf{break}
			\EndIf
			\State $z_{\text{prev}} \leftarrow z^{(t)}$ 
			\ForAll{edges $(i,j) \in E$ \textbf{in parallel}}
				\State For each commodity $k \in \mathcal{K}$, $\phi_{ij,k} \leftarrow h_{ik}^{(t)} - h_{jk}^{(t)} - \psi_{ij}^{(t)}$.
				\State For each commodity $k \in \mathcal{K}$, $\Delta f_{ij,k} \leftarrow \max(-f_{ij,k}^{(t)}, \beta_{ij} \cdot \phi_{ij,k})$.
				\State Estimate objective change $\Delta z_{ij}$ if update $\{\Delta f_{ij,k}\}$ is applied.
				
				\If{$\Delta z_{ij} \le 0$} 
					\State For each $k \in \mathcal{K}$, \textcolor{blue}{$v_{ij,k}^{(t+1)} \leftarrow \gamma \cdot v_{ij,k}^{(t)} + \Delta f_{ij,k}$}.
					\State For each commodity $k \in \mathcal{K}$, \textcolor{blue}{$f_{ij,k}^{(t+1)} \leftarrow \max(0, f_{ij,k}^{(t)} + v_{ij,k}^{(t+1)})$}.
				\Else
					\State For each commodity $k \in \mathcal{K}$, $f_{ij,k}^{(t+1)} \leftarrow f_{ij,k}^{(t)}$, \textcolor{blue}{$v_{ij,k}^{(t+1)} \leftarrow v_{ij,k}^{(t)}$}.
					\State $\beta_{ij} \leftarrow \max(\beta_{min}, \beta_{ij} \cdot \alpha_{dec})$.
				\EndIf
				\If{$t \pmod{10} = 0$} 
					\State $\beta_{ij} \leftarrow \beta_{ij} \cdot \alpha_{inc}$ 
				\EndIf
			\EndFor
			
			\State $t \leftarrow t + 1$
		\EndWhile
		\State \textbf{Return} final flow $f^{(t)}$
		\end{algorithmic}
	\end{algorithm}

	\section{Further Discussion on Optimality Condition}
	As commonly understood, directed edges correspond to pipes of the physical network and vertices  are pipe junctions. In practice, for any physical network, it is impossible for the flow in a pipe to exceed its capacity . In this section we would discuss the optimality condition for 
	multicommodity flow problem when only relaxing the flow conservation constraints.

	The basic formulation  is as follows:
	\begin{equation}\label{feq5}
	\begin{aligned}
	{\bf min}\quad &z =  \sum_{i,k}\int_{0}^{\Delta_{ik}}h_{ik}(\omega)d\omega \\
	{\bf s.t}\quad & \sum_{ k\in \mathcal{K}}f_{ij,k}\leq u_{ij}, \forall (i,j) \in E\\
	 \quad& f_{ij,k}\geq 0, \forall k\in \mathcal{K}, (i,j) \in E\\
	\end{aligned}
	\end{equation}
	where \(h_{ik}\) is the height function and \(\Delta_{ik}\) is defined as Expression~(\ref{eq3}).

	In the programming above, the objective function is the sum of the integrals of the vertex height functions. The   flow conservation constraints are relaxed but the  capacity constraints are maintained. Obviously, the feasible region of Expression (\ref{eq1}) is  not empty if and only if the minimum value of the objective function of Programming~(\ref{feq5}) is zero; the feasible region of Expression (\ref{eq1}) is   empty if and only if the minimum value of the objective function of Programming~(\ref{feq5}) is greater than zero. Since  Programming~(\ref{feq5}) is convex,  the solution \(\boldsymbol{f}^\ast\) of Programming~(\ref{feq5}) is optimal if and only if it satisfies the  following Karush-Kuhn-Tucker conditions:

	\begin{equation}\label{eq61}
	\begin{aligned}
	&{\bf Stationarity}\\
	& \quad -\frac{\partial z}{\partial f_{ij,k}} = -\mu_{ij,k} +  \gamma_{ij},\forall k\in \mathcal{K}, (i,j) \in E\\
	&{\bf Primal~feasibility}\\
	& \quad \sum_{ k\in \mathcal{K}}f_{ij,k}\leq u_{ij}, \forall (i,j) \in E\\
	& \quad -f_{ij,k}\leq 0, \forall k\in \mathcal{K}, (i,j) \in E\\
	&{\bf Dual~feasibility}\\
	& \quad \mu_{ij,k} \geq 0,\forall k\in \mathcal{K}, (i,j) \in E\\
	& \quad \gamma_{ij} \geq 0, \forall (i,j) \in E\\
	&{\bf Complementary~slackness}\\
	& \quad \mu_{ij,k}f_{ij,k} = 0,\forall k\in \mathcal{K}, (i,j) \in E\\
	& \quad \gamma_{ij}(\sum_{ k\in \mathcal{K}}f_{ij,k} - u_{ij}) = 0,\forall (i,j) \in E.\\
	\end{aligned}
	\end{equation}
	
	Obviously,
	\begin{equation}
	\label{cpgrad1}
	\begin{aligned}
	\frac{\partial z}{\partial f_{ij,k}}  &=  h_{ik}(\Delta_{ik}) \frac{\partial \Delta_{ik}}{\partial f_{ij,k}} + h_{jk}(\Delta_{jk}) \frac{\partial \Delta_{jk}}{\partial f_{ij,k}}\\
		&= (-h_{ik}(\Delta_{ik})) + (h_{jk}(\Delta_{jk})) \\
		&=  h_{jk} - h_{ik}.\\
	\end{aligned}
	\end{equation}
	Substituting the  expression above into Stationarity expression in KKT conditions,
	\begin{equation*}
	  \gamma_{ij} = \mu_{ij,k}  + h_{ik} - h_{jk},\forall k\in \mathcal{K}, (i,j) \in E.\\
	\end{equation*}
	
	If edge\((i,j)\) is not saturated, by complementary slackness, \( \gamma_{ij} = 0\). If edge\((i,j)\) is  saturated,  \( \gamma_{ij} \geq  0\). Intuitively, \( \gamma_{ij} \) could  be interpreted as the pipe pressure on edge\((i,j)\). If edge\((i,j)\) is not saturated, the \emph{pipe pressure} is zero; otherwise it is greater than or equal to zero.

	For any used edge of commodity \(k\), i.e. \(f_{ij,k} > 0\), by complementary slackness \( \mu_{ij,k}f_{ij,k} = 0\) in KKT conditions, we have \( \mu_{ij,k} = 0\). Therefore,
	\begin{equation*}
	  \gamma_{ij} =  h_{ik} -h_{jk},\forall k\in \mathcal{K}, (i,j) \in E.\\
	\end{equation*}
	That is, the height difference between vertex \(i\) and vertex \(j\) for commodity \(k\) is equal to the pipe pressure \( \gamma_{ij}\) of edge \((i,j)\).
	
	That is, for any commodity \(k\) using edge \((i,j)\), its height difference must equal the pipe pressure  \( \gamma_{ij}\). A critical consequence of this result is that \emph{for any saturated edge, the height difference \(h_{ik} - h_{jk}\) must be the same for all commodities \(k\) that have a positive flow on that edge.}

	For any unused edge of commodity \(k\), due to \(\mu_{ij,k} \geq 0\), we have 
	\begin{equation*}
	\begin{aligned}
	  \gamma_{ij} &=  \mu_{ij,k} + h_{ik} -h_{jk} \geq  h_{ik} -h_{jk}, \\
	  &\forall k\in \mathcal{K}, (i,j) \in E\\
	\end{aligned}
	\end{equation*}
	
	That is, the height difference between vertex \(i\) and vertex \(j\) for commodity \(k\) is less than or equal to the pipe pressure \( \gamma_{ij}\).

	Note that the pipe pressure \( \gamma_{ij}\) is zero when the edge is unsaturated. Therefore, the optimality conditions can be concluded as follows:
	\begin{enumerate}[label=(\roman*)]
	\item for any used  edge of commodity \(k\), if it is unsaturated, the height difference between vertex \(i\) and vertex \(j\) for commodity \(k\) is zero;
	\item for any used edge of commodity \(k\), if it is saturated, the height difference between vertex \(i\) and vertex \(j\) for commodity \(k\) is equal to the pipe pressure \( \gamma_{ij}\);
	\item for any unused edge of commodity \(k\), the height difference between vertex \(i\) and vertex \(j\) for commodity \(k\) is less than or equal to the pipe pressure \( \gamma_{ij}\).
	\end{enumerate}

	\section{Numerical Results}
	
	In this section, we empirically evaluate the performance of the proposed algorithms on a suite of standard benchmark problems. We first describe the test instances and then present a comparative analysis of the results.
	
	\subsection{Test Problems}
	
	Our numerical experiments were conducted on three sets of well-known benchmark problems for the multicommodity flow problem\citep{babonneau2006solving,larsson2004augmented}.
	
	\begin{enumerate}
		\item \textbf{Planar Problems:} This set, generated by \cite{larsson2004augmented}, comprises 10 instances designed to simulate telecommunication networks. Vertices are placed randomly in a plane, and arcs connect neighboring vertices to form a planar graph. Commodities, demands, and capacities are chosen randomly from uniform distributions.
		\item \textbf{Grid Problems:} This set contains 15 instances with a grid network structure, where most vertices have a degree of four. This topology is characterized by a large number of alternative paths between vertex pairs.
		\item \textbf{Telecommunication Problems:} This collection includes real-world and realistic problems of varying scales. Notably, it features the \texttt{904} instance, a large-scale network with 904 arcs and 11,130 commodities, as referenced in the survey by \cite{ouorou2000survey}.
	\end{enumerate}
	
	These problem sets are publicly available\footnote{\url{http://www.di.unipi.it/di/groups/optimize/Data/MMCF.html}} and standard for benchmarking multicommodity flow algorithms. Since the objective of our study is to find a \emph{feasible} flow (i.e., a zero-stable pseudo-flow), the arc cost data present in the original benchmarks were \emph{disregarded}. For each problem, we report the number of vertices $|V|$, arcs $|E|$, and commodities $|\mathcal{K}|$.
	
	\subsection{Numerical Experiments}
	
	We implemented and compared the three primary algorithmic variants discussed in this paper:
	\begin{itemize}
		\item \textbf{\texttt{PDR-ESO} (PDR via Exact Subproblem Optimization):} The method from Algorithm \ref{alg:qp_projection}, which solves the edge subproblem (Expression (\ref{eq9})) exactly in each iteration.
		\item \textbf{\texttt{PDR-AGD} (PDR via Adaptive Gradient Descent):} The heuristic method from Algorithm \ref{alg:gradient_descent}, which uses a projected gradient step with an adaptive learning rate.
		\item \textbf{\texttt{PDR-GDM} (PDR via Gradient Descent with Momentum):} The enhanced heuristic from Algorithm \ref{alg:momentum_gradient}, which incorporates a momentum term to accelerate and stabilize convergence.
	\end{itemize}
	
	\subsubsection{Experimental Setup}
	All algorithms were implemented in CUDA and executed on a NVIDIA A100 GPU. The algorithms were tested on 28 benchmark instances, with performance primarily evaluated based on the number of iterations and total computation time required to satisfy the convergence criterion $|z^{(t)} - z^{(t-1)}| < 10^{-6}$. Here, $z^{(t)}$ is the objective value (from Expression (\ref{eq8}), representing the sum of squared constraint violations) at iteration $t$.
	
	For \textbf{\texttt{PDR-AGD}}, the parameters were configured as follows: the initial learning rate $\beta_{\text{init}}$ was set to 0.25, with a minimum learning rate $\beta_{\text{min}} = 1.0 / |\mathcal{K}|$ and a maximum $\beta_{\text{max}} = 0.5$. The learning rate adjustment factors were $\alpha_{\text{inc}} = 2.0$ and $\alpha_{\text{dec}} = 0.5$.
	
	For \textbf{\texttt{PDR-GDM}}, the parameters were set with an initial learning rate $\beta_{\text{init}} = 0.1$, momentum coefficient $\gamma = 0.9$, a minimum $\beta_{\text{min}} = 1.0 / |\mathcal{K}|$, and a maximum $\beta_{\text{max}} = 0.2$. The adjustment factors $\alpha_{\text{inc}}$ and $\alpha_{\text{dec}}$ remained 2.0 and 0.5, respectively. The complete experimental results are presented in Table~1.
		
		\begin{table}[h!]
		\centering
		\caption{Complete performance comparison of \texttt{PDR-ESO}, \texttt{PDR-AGD}, and \texttt{PDR-GDM} algorithms.}
		\label{tab:full_results}
		\resizebox{\textwidth}{!}{%
		\begin{tabular}{@{}l
						S[table-format=5.0]
						S[table-format=5.0]
						S[table-format=5.0]
						@{\qquad} 
						S[table-format=5.0]
						S[table-format=4.4]
						S[table-format=1.7]
						@{\qquad} 
						S[table-format=5.0]
						S[table-format=4.4]
						S[table-format=1.7]
						@{\qquad} 
						S[table-format=4.0]
						S[table-format=3.4]
						S[table-format=1.7]
						@{}}
		\toprule
		& & & & \multicolumn{3}{c}{\textbf{PDR-ESO  Algorithm}} & \multicolumn{3}{c}{\textbf{PDR-AGD Algorithm}} & \multicolumn{3}{c}{\textbf{PDR-GDM Algorithm}} \\
		\cmidrule(lr){5-7} \cmidrule(lr){8-10} \cmidrule(lr){11-13}
		\textbf{Instance} & {\textbf{$|E|$}} & {\textbf{$|V|$}} & {\textbf{$|\mathcal{K}|$}} & {\textbf{Iter.}} & {\textbf{Time (s)}} & {\textbf{Obj $z$}} & {\textbf{Iter.}} & {\textbf{Time (s)}} & {\textbf{Obj $z$}} & {\textbf{Iter.}} & {\textbf{Time (s)}} & {\textbf{Obj $z$}} \\
		\midrule
		C22 & 22 & 14 & 23 & 56 & 0.0064 & 0.0000023 & 61 & 0.0053 & 0.0000028 & 115 & 0.0102 & 0.0000195 \\
		C148 & 148 & 61 & 122 & 1075 & 0.1451 & 0.0001331 & 1734 & 0.1668 & 0.0010229 & 399 & 0.0389 & 0.0014546 \\
		C904 & 904 & 106 & 11107 & 120 & 0.1476 & 0.0000061 & 126 & 0.0670 & 0.0000058 & 175 & 0.0959 & 0.0000208 \\
		Cgd1 & 80 & 25 & 50 & 374 & 0.0458 & 0.0000203 & 441 & 0.0418 & 0.0000245 & 176 & 0.0190 & 0.0000604 \\
		Cgd2 & 80 & 25 & 100 & 1915 & 0.2586 & 0.0001648 & 7282 & 0.6536 & 0.0012701 & 212 & 0.0193 & 0.0000083 \\
		Cgd3 & 360 & 100 & 50 & 1143 & 0.1595 & 0.0000685 & 1168 & 0.1145 & 0.0000705 & 190 & 0.0187 & 0.0000169 \\
		Cgd4 & 360 & 100 & 100 & 2108 & 0.3106 & 0.0001692 & 2615 & 0.2521 & 0.0002245 & 820 & 0.0798 & 0.0001075 \\
		Cgd5 & 840 & 225 & 100 & 2523 & 0.3990 & 0.0001611 & 2536 & 0.2556 & 0.0001633 & 622 & 0.0658 & 0.0000365 \\
		Cgd6 & 840 & 225 & 200 & 3484 & 0.6091 & 0.0003996 & 5295 & 0.5536 & 0.0005358 & 1391 & 0.1442 & 0.0037313 \\
		Cgd7 & 1520 & 400 & 400 & 4595 & 1.8731 & 0.0002962 & 4610 & 0.5137 & 0.0002983 & 1243 & 0.1448 & 0.0000896 \\
		Cgd8 & 2400 & 625 & 500 & 7273 & 4.1351 & 0.0005089 & 7893 & 0.9848 & 0.0006005 & 1814 & 0.2363 & 0.0050351 \\
		Cgd9 & 2400 & 625 & 1000 & 6481 & 4.1993 & 0.0169427 & 12614 & 1.9144 & 0.0013713 & 2534 & 0.4277 & 0.0023245 \\
		Cgd10 & 2400 & 625 & 2000 & 7820 & 6.2245 & 0.0006290 & 8742 & 2.2619 & 0.0006594 & 2453 & 0.6809 & 0.0099367 \\
		Cgd11 & 2400 & 625 & 4000 & 7954 & 6.8848 & 0.0004862 & 7963 & 3.4184 & 0.0004857 & 2057 & 0.9748 & 0.0001125 \\
		Cgd12 & 3480 & 900 & 6000 & 11318 & 34.0878 & 0.0007265 & 11338 & 9.0589 & 0.0007165 & 2977 & 2.6477 & 0.0001670 \\
		Cgd13 & 3480 & 900 & 12000 & 11803 & 55.6226 & 0.0007470 & 11822 & 16.8680 & 0.0007317 & 3092 & 5.0515 & 0.0001691 \\
		Cgd14 & 4760 & 1225 & 16000 & 15864 & 118.0360 & 0.0010348 & 15872 & 39.6654 & 0.0010232 & 4192 & 11.9781 & 0.0002360 \\
		Cgd15 & 4760 & 1225 & 32000 & 16470 & 240.2010 & 0.0011822 & 16520 & 80.6890 & 0.0010956 & 4352 & 24.4519 & 0.0002394 \\
		Cpl30 & 150 & 30 & 92 & 340 & 0.0473 & 0.0000229 & 359 & 0.0325 & 0.0000244 & 140 & 0.0132 & 0.0001120 \\
		Cpl50 & 250 & 50 & 267 & 669 & 0.1315 & 0.0000496 & 880 & 0.0856 & 0.0000790 & 173 & 0.0174 & 0.0000094 \\
		Cpl80 & 440 & 80 & 543 & 620 & 0.1767 & 0.0003072 & 1040 & 0.1092 & 0.0001094 & 202 & 0.0212 & 0.0000133 \\
		Cpl100 & 532 & 100 & 1085 & 1048 & 0.3735 & 0.0001003 & 1362 & 0.1541 & 0.0001146 & 239 & 0.0278 & 0.0000095 \\
		Cpl150 & 850 & 150 & 2239 & 1735 & 1.1987 & 0.0055658 & 37873 & 4.9982 & 0.0287100 & 1358 & 0.1977 & 0.0000120 \\
		Cpl300 & 1680 & 300 & 3584 & 2719 & 2.1256 & 0.0002071 & 2733 & 0.7969 & 0.0002074 & 682 & 0.2213 & 0.0000447 \\
		Cpl500 & 2842 & 500 & 3525 & 3222 & 3.0401 & 0.0002559 & 3228 & 1.4272 & 0.0002559 & 829 & 0.4028 & 0.0000576 \\
		Cpl800 & 4388 & 800 & 12756 & 5345 & 35.7495 & 0.0004423 & 5473 & 10.7503 & 0.0004566 & 1358 & 3.0862 & 0.0001911 \\
		Cpl1000 & 5200 & 1000 & 20026 & 9047 & 109.7860 & 0.0037456 & 10695 & 38.5915 & 0.0008608 & 2298 & 9.4954 & 0.0020320 \\
		Cpl2500 & 12990 & 2500 & 81430 & 15225 & 1845.3300 & 0.0211950 & 21408 & 788.5840 & 0.0018430 & 4549 & 191.1520 & 0.0050321 \\
		\bottomrule
		\end{tabular}
		} 
		\end{table}

	\subsubsection{Analysis of Convergence Behavior}
	
	The results clearly demonstrate the superior performance of the momentum-accelerated method. The \texttt{PDR-GDM} algorithm consistently converges in the fewest iterations and, consequently, the shortest time across the vast majority of instances. As a representative example, on instance \texttt{Cpl1000}, \texttt{PDR-GDM} converges in 2298 iterations, whereas the standard \texttt{PDR-AGD} method requires 10695 iterations---a 4.6-fold improvement.

	An insightful and arguably counter-intuitive result is that the \texttt{PDR-ESO} algorithm, despite solving each edge subproblem to optimality, frequently requires a substantially greater number of iterations to converge than its approximation-based counterpart, \texttt{PDR-GDM}. This observation underscores the limitations of a purely greedy optimization approach. It leads to the pivotal conclusion that local optimality at each step does not guarantee, and may not even facilitate, rapid global convergence. Consequently, a sequence of less precise yet more aggressive steps, such as those generated by \texttt{PDR-GDM}, can be demonstrably more efficient, ultimately achieving convergence with a lower total iteration count.
	
	\subsubsection{Computational Time and Solution Quality}

	In terms of computational time, \textbf{\texttt{PDR-GDM} consistently emerges as the most efficient algorithm.} This superior performance is a direct result of its ability to significantly reduce the number of iterations while maintaining a computationally inexpensive update step. In contrast, \textbf{\texttt{PDR-ESO} is almost always the slowest}, hampered by the high computational cost of its exact solver ($O(|\mathcal{K}|\log |\mathcal{K}|)$ per edge). While \texttt{PDR-AGD} is generally faster than \texttt{PDR-ESO}, it cannot match the accelerated convergence of the momentum-based approach.
	
	The practical impact of this efficiency is particularly evident on large-scale instances. For example, the \textbf{\texttt{Cgd15} problem} involves \textbf{over 152 million} variables ($|E| \times |\mathcal{K}| = 4760 \times 32000$). Remarkably, \texttt{PDR-GDM} finds a feasible solution for this massive instance in just \textbf{24.45 seconds}.
	
	Regarding solution quality, all three algorithms successfully reduce the objective value $z$---representing the sum of squared constraint violations---to a very small magnitude. It is important to contextualize this final objective value. Given the enormous number of variables and constraints in these benchmark problems, the final values of $z$ indicate that the \textbf{ violation per single constraint} (both for arc capacities and flow conservation) is \textbf{negligibly small}. Therefore, for all practical purposes, the solutions can be considered \textbf{fully feasible and are of high quality}, satisfying the primary goal of the study.

	\section{Conclusion}
	
	In this paper, we introduced a novel theoretical framework for the multicommodity flow problem, grounded in the physical intuition of potential (height) and congestion. By defining vertex-specific potential functions and edge-specific congestion functions, we established a simple and intuitive set of optimality conditions. We then reformulated the classical feasibility problem as a convex optimization problem with nonnegativity constraints, where the existence of a feasible solution is equivalent to the objective function attaining a minimum value of zero.
	
	A key advantage of this framework is its edge-separable structure, which makes it particularly well-suited for parallel and distributed computation. Leveraging this property, we proposed three algorithms: a method that solves an exact quadratic subproblem for each edge, and two computationally efficient heuristics based on adaptive gradient descent and its momentum-accelerated variant. Numerical experiments on benchmark instances demonstrated that our algorithms are highly efficient for large-scale networks, with the momentum-based heuristic (\texttt{PDR-GDM}) exhibiting particularly rapid convergence.
	
	While the proposed framework and algorithms show great promise for solving large-scale problems, this work also opens several avenues for future research. The most pressing theoretical question is the formal analysis of the worst-case convergence rates and computational complexities of our proposed algorithms. Furthermore, exploring more sophisticated adaptive learning rate and momentum schemes could lead to even faster convergence. Finally, applying this potential and congestion framework to other classes of network flow problems, such as those with nonlinear cost functions or additional side constraints, presents a rich area for further investigation.

	
	\section{Code and Data Disclosure}\label{sec:Code and Data Disclosure}The code and data to support the numerical experiments in this paper can be found at \url{https://github.com/mrgump-123/Localized-Multicommodity-Flow}.

\bibliographystyle{plainnat}
\bibliography{references}

\end{document}